\documentclass{amse-new}

\graphicspath{{figures/}}

\numberwithin{equation}{section} 

\begin{document}

 \PageNum{1}
 \Volume{201x}{Sep.}{x}{x}
 \OnlineTime{August 15, 201x}
 \DOI{0000000000000000}
 \EditorNote{Received x x, 201x, accepted x x, 201x}

\abovedisplayskip 6pt plus 2pt minus 2pt \belowdisplayskip 6pt
plus 2pt minus 2pt
\def\vsp{\vspace{1mm}}
\def\th#1{\vspace{1mm}\noindent{\bf #1}\quad}
\def\proof{\vspace{1mm}\noindent{\it Proof}\quad}
\def\no{\nonumber}
\newenvironment{prof}[1][Proof]{\noindent\textit{#1}\quad }
{\hfill $\Box$\vspace{0.7mm}}
\def\q{\quad} \def\qq{\qquad}
\allowdisplaybreaks[4]


\AuthorMark{Minzheng Li }                             

\TitleMark{On the Equivalence of the Newton-Raphson Algorithm and PDE of Conservation of Electric Charge}  

\title{On the Equivalence of the Newton-Raphson Algorithm and PDE of Conservation of Electric Charge        
\footnote{Supported by \ldots (Grant No. \ldots)}}                  

\author{Minzheng Li}             
    {Address School of Mathematical Sciences, Peking University \\
    E-mail 1801110057@pku.edu.cn \,$:$ }

\maketitle%

\Abstract{The main result characterises the equivalence of the Newton-Raphson algorithm and PDE of conservation of electric charge. Based on this equivalence we analyse the properties of the Fisher-scoring method, a variant of Newton's method commonly used in statistics. Conjecture of embedding physics models to solve Newton's method is proposed.}      


\MRSubClass{68Wxx,00A79}      

\section{Introduction}
The Newton-Raphson algorithm, or in the popular name, Newton's method, is a primary and ancient method of numerically minimisation or identifying roots of a target function. Given a function $ F $ with its gradient matrix $ J $:$$
J = \mathrm{grad}(F)
$$Newton suggests the following iteration:$$
x_{k+1}=x_k - J^{-1}(x_k)F(x_k)
$$in the continuous-time form this is written as:$$
\dot{x} = - J^{-1}(x)F(x)
$$The method is inspired by a simple Taylor expansion: $$
0= F(x^*) = F(x) + J(x)(x^*-x)
$$ where $ x^* $ is the root of $ F $ and $ x $ is near $ x^* $, and remains one of the most basic numerical tools over hundreds of years. However, as mathematical models and statistical models retrieving more and more complexity, more challenges are facing this traditional numerical algorithm. Typically in theory the sufficient and necessary condition of its convergence remain obscure to mathematical researches.

In practice the effectiveness of Newton's method is struggled with the need of locating a good start point and the necessity of the target function being not too rough. In coping with the start point issue, Fisher-scoring method, a variant of Newton's method, simply selects the ordinary Least Square estimator as its start point. In machine learning practice this is usually tried by equilibriously choosing many trial start points. In coping with the roughness, typically in many mathematical and physical models, line searches and inexact solvers may be provided\cite{LI2022107296}\cite{doi:10.1137/15M1028078}\cite{LAMPRON2021114091}\cite{YANG20171}. There are also many other techniques in applying Newton's method in real problems. Fisher-scoring algorithm change Jacobi matrix into Fisher matrix. \cite{MORADIAN2021109886} and \cite{SEREETER2019157} also succeed in applying variant substitutes of Jacobi matrices according to the nature of the problems. Alternating minimization scheme can apply\cite{LAMPRON2021114091}. Last but not the least, Newton-GMRes methods\cite{doi:10.1137/S1064827599363976}\cite{articleGMRes} revert determining step-lengths to solving a reverted algebraic equation.

\cite{generalisedNR} advances theoretically by suggesting considering iterations of the form$$
x^{(k+1)} = s^{-1}\left( s(x^{(k)}) - J_s(x^{(k)})[J_f(x^{(k)}]^{-1} f(x^{(k)}) \right)
$$which is equivalent to optimise $ F = f \circ s^{-1} $, and thus the convergence region can be expanded. However it remains to be found out what $ s(\cdot) $ function to be used. Taking $ s(\cdot) $ for example a constant multiplying identity matrix, the method is equivalent to multiplying a constant in every iteration step-length. We may brutally say, their text results in manipulating step-lengths in practice, and as will be shown that ours results in manipulating differentiation grid length compared.

The equivalence of Newton's method with the PDE of conservation of electric charge is proposed and proved in this text. First we prove a lemma in statistical physics. We have questioned the origination of this result. Fluid dynamics and statistical physics both have hundreds of years of history, but after enquiring physics texts we fail to find similar proof. We believe the work might have been done but discarded as trivial. We showcase the main result equivalence by the lemma and some additional probabilistic analysis. Albeit with much higher computational costs, solving PDE observes an overall picture of the behaviour of Newton's method, and thus free us from the necessity of locating good start points. Favoured as the intriguing and profound equivalent form, this mathematical discovery pledges to be more deeply studied by expertise from various backgrounding. In the end of this paper, total momenta of Fisher-scoring method is defined, though it is far from proving its convergence. Also Newton's method is conjectured to be embedded into many physics models.

By proposing this equivalence, we in a large degree revert the convergence issue of the Newton's method. Consider removing the Jacobi matrix and the iteration:$$
\dot{x} = - F(x)
$$By running equivalent PDE, as long as the numerical lattice is sufficiently small, we capture every root of $ F $ where $ J $ is locally convex, and it is a case that we prove for the Fisher-scoring algorithm. However it is a trade-off that how small the lattice we should choose. We tackle the roughness of the function provided and tend to choose smaller lattice, but it takes costs. On the one hand, smaller lattice requires higher computational cost. Consider three dimension which is low dimensional scenario mathematically. By choosing 100 grids in each dimension we take $ 5^3 = 125 $ times computation than choosing 20 grids each. If further consider four dimension this becomes $ 5^4 = 625 $ times, an unbearable additional burden of computation. On the other hand, smaller lattice may suffer from round-off error. If the grid length is $ 0.1 $ we do not really care about numerical round-off, but if grid length $ 0.001 $ it is hardly say that round-off effect is not an issue.

Our research starts from image processing. For decades mathematicians and computer scientists have exhausted imagination to explore techniques, but image processing remain unsatisfactory with regard to requirements of computer vision. We started by developing the PDE of conservation of electric charge $$
\frac{\partial\rho}{\partial t} + \mathrm{div}\mathbf{j}=0
$$into some image processing techniques, and ended finding the relationship of particles moving in the velocity field and its statistical description.

Macro physics and particle physics are two systems in the study. Electromagnetic theories and fluid dynamic theories study macro statistical phenomena, and classical dynamics and quantum mechanics study the existence and movement of particles. Fluid dynamics do not consider non-physic scenario like 'a small iron ball flowing in the air stream'. The so-called continuity equation in fluid dynamics is based on the continuity assumption, that the media is continuous and uniformly shaped. And in classic Navier-Stokes equations, the velocity functions are not real velocity of particles but statistical description of the macro. However in this text, 'A small iron ball flowing in air stream' is what we do consider, violating physics but meaning in mathematics.

On the other hand there are many results in physics linking the micro with the macro. The Liouville theorem in statistical physics targets on Hamilton-conserving particles. 
Boltzmann equations, usually written as $$
\frac{\partial}{\partial t}\rho(x,v,t) + v\cdot \frac{\partial}{\partial x}\rho(x,v,t) + F\cdot \frac{\partial}{\partial v}\rho(x,v,t) = -\tau (\rho -\rho^{(eq)})
$$ $ \rho^{(eq)} $ being local equilibrium state, appears in many modelling problems. DFT(Density Functional Theory)-based quantum chemistry enables scientists to simulate the structure of crystals or the molecular dynamics of liquids from the quantum level. \cite{PhysRev.140.A1133} proposed the well-known Kohn-Sham equations, regarded as a cornerstone of DFT but taking too high computational costs. \cite{osti_1542041} takes a different approach by learning potential functions into neural networks, which is a revolution in the study of First Principles in computational physics. By the time the text is written no connection with these previous results has been found.

The outline of this paper is as follows: In the section 'A Lemma in Statistical Physics' a lemma is half-proved. By saying half-proved, the particle nature containing statistical nature is proved, and it is used directly in the main result. However the statistical nature containing particle nature is non-strictly illustrated, leaving the work of probabilistic analysis in Main Result section. A simple GLM solving MLE numerical experiment is carried out next. Finally many related studies are given in the last section 'Extensions'.

Minimisation and root solving are basic mathematical problems. Computational mathematics have developed rich methods from the time modern computers were invented. Most parallel extensions focus on solving specific problems. Advance in techniques confuses us more on the essence of science, and it is believed that only by discovering the essence of science can mathematicians reveal true spirit of real techniques. Considering particle-statistics equivalence is a small step forward theoretically. It is far from being as important as say like Fourier transform that equals a measurable function into completely Fourier spaces, but we do hope our particle-statistics equivalence gives us different perspectives on the phenomenon of the Newton's method, and possibly advances new techniques.

\section{Notation convention}
In the whole text $ \mathrm{div} $ and $ \triangledown $ denote taking divergence and taking gradient of a function or a vector function of its location part $ x $ but not including its time part $ t $.

\section{Preliminary}
\subsection{The PDE of Conservation of Electric Charge}
Conservation of electric charge is the basic principle of the electromagnetics, with its equations being the precondition for Maxwell equations. Let $ \rho(x,t) $ be the electric density of location-time $ (x,t) $, defined as for any location domain $ \Omega $, $ \int_\Omega \rho dx $ is the total electric charge within. Let $ \mathbf{j}(x,t) $ be electric current, defined as $ \int_{[t_1,t_2]}\int_S \mathbf{j}\cdot\mathbf{n}~dSdt $ being the electric charge passed through the surface $ S $ within time $ [t_1,t_2] $. The conservation of electric charge states: \textbf{electric charges do not appear or disappear in vain}. Written in mathematics by a simple Gauss-Green theorem, it is
$$
\frac{\partial\rho}{\partial t} + \mathrm{div}\mathbf{j}=0
$$

Generalised in fluid dynamics, the PDE of conservation of electric charge appears as the continuity equations. Revert the electric current $ \mathbf{j} $ as the product of the density $ \rho(x,t) $ and the velocity $ \mathbf{v}(x) $. After the same Gauss-Green theorem the continuity equation writes as

\begin{equation}\label{continuity_equation}
\frac{\partial\rho(x,t)}{\partial t} + \mathrm{div}(\rho(x,t)\mathbf{v}(x))=0
\end{equation}
To give an example of \eqref{continuity_equation}, consider clouds in the sky. $ \rho $ is the density of clouds, $ \mathbf{v} $ is statistical velocity, $ x $ is two dimensional. The readers might wonder whether this can be written as $ \frac{\partial\rho(x,t)}{\partial t} + \mathrm{div}(\mathbf{v}(x))\rho(x,t)=0 $, kicking off $ \rho $ from the divergence operator. In physics this is no go. In the numerical study later this is also a no go.

Notice here the velocity $ \mathbf{v}(x) $ is statistical. In physics we mean statistically we observe this velocity in the macro, which is different than later. We will prove in the whole text that if a mathematical phenomenon is scripted by a recipe like this, then continuity equations designates a true velocity of particle. In physics we consider a continuity equation to describe physical worlds, but here we consider a statement given by the ensemble of continuity equations.

It is worth to mention that the phenomenon characterised by the PDE of conservation of electric charge or the continuity equation is everywhere in physics. In quantum mechanics probability current density theory is given by \begin{equation}\label{Schrodinger1}
\frac{\partial}{\partial t}\rho + \frac{\hbar}{i}\frac{1}{2m}\mathrm{div}\left[  \psi^*\triangledown\psi- \psi\triangledown \psi^* \right] = 0
\end{equation}
which is a simple corollary from the Schr\"{o}dinger equation.

\subsection{The Fisher-Scoring Algorithm}
This subsection is a brief review of mathematical statistics. Readers uninterested with mathematical statistics can skip.

In the most general statistics notation, let $ f $ be probability density function, $ l(\beta) =\log f $ be log-likelihood function, and the partial derivative of which, $ S(\beta) = \frac{\partial l}{\partial \beta} $, is named as the score function. $ S(\beta)=0 $ yields a solution maximum likelihood estimator.

By taking derivatives of $ \beta $ at both sides of $ \int f dx=1 $,
\begin{lemma}
(The first Bartlett identity) the expectation of Score with true $ \beta $ substituted in is 0 $$
ES(\beta)=0
$$
\end{lemma}

By further taking derivatives at both sides of the first Bartlett identity,
\begin{lemma}\label{Bartlett2}
(The second Bartlett identity) $$
E\left(-\frac{\partial S}{\partial \beta^T}\right)= I(\beta)
$$
\end{lemma}

The GEE theory regard a wealth of statistical modelling to take in forms of $$
S(\beta) =\sum\limits_{j=1}^n \psi_j(X_j,\beta)
$$from which the consistency is of independent and non-identically distributed theorems of large number, and the asymptotic normality is of the Lindeberg-Feller theorem. See \cite[Section 5.4]{Shaojun}. The classic Newton-Raphson algorithm \begin{equation}\label{Newton_Raphson}
\widehat{\beta}^{(l+1)}=\widehat{\beta}^{(l)}-\left(\frac{\partial S}{\partial \beta^T}\right)^{-1}(\widehat{\beta}^{(l)})S(\widehat{\beta}^{(l)})
\end{equation} does work for solving maximum likelihood estimation, but what is usually preferred is the so-called Fisher-scoring algorithm \begin{equation}\label{Fisher_Scoring}
\widehat{\beta}^{(l+1)}=\widehat{\beta}^{(l)}+ I(\widehat{\beta}^{(l)})^{-1}S(\widehat{\beta}^{(l)})
\end{equation}which substitutes the derivative of the score by its expectation with negative charge, namely Fisher information matrix. We grace this variant of Newton's method as it embraces better analytic properties shown in the ending section.

As a special example, by taking particularly canonical link in generalised linear models, the Fisher-scoring algorithm is equivalent to Newton's method. In the most general statistical notations $ Y\sim \exp(\frac{y\theta - b(\theta)}{a(\phi)}) $, $ \theta = u(x^T\beta) $.
$$
likelihood = \exp\left( \sum\limits_{j=1}^n \frac{y_j\theta_j - b(\theta_j)}{a(\phi_j)} \right)
$$ $$
loglikelihood = \sum\limits_{j=1}^n \frac{y_j\theta_j - b(\theta_j)}{a(\phi_j)}
$$ $$
Score = \frac{\partial loglikelihood}{\partial\beta} = \sum\limits_{j=1}^n \frac{y_j - \dot{b}(\theta_j)}{a(\phi_j)}\dot{u}(x_j^T\beta)x_j
$$ $$
H = \frac{\partial Score}{\partial\beta^T} = \sum\limits_{j=1}^n x_j^T \dot{u}(x_j^T\beta)^T \frac{(- \ddot{b}(\theta_j))}{a(\phi_j)}\dot{u}(x_j^T\beta)x_j + J
$$ where $ J $ is the item taking derivatives of $ \beta $ from $ \dot{u}(x_j^T\beta) $, satisfying $ EJ =0$. Generally we have$$
I(\beta) = E(-H(\beta)) =\sum\limits_{j=1}^n x_j^T \dot{u}(x_j^T\beta)^T \frac{(\ddot{b}(\theta_j))}{a(\phi_j)}\dot{u}(x_j^T\beta)x_j
$$ and in cases of canonical links $ J=0 $ and thus $ I = -H $.

\section{A Lemma in Statistical Physics}
Consider particles moving in the velocity field $ \mathbf{v}(\cdot) $ with the corresponding ODE $$
\frac{d}{dt}x = \mathbf{v}(x)
$$ the existence and uniqueness of solutions of which is characterized by the classic Pichard iteration theorem. Consider the two statements:

\begin{itemize}
\item Statement One (particle description of fluid): Particles $ \widehat{\beta} $'s move in velocity field $ \mathbf{v}(x) $. Written in differentiation this is $$\widehat{\beta}^{(t+dt)}=\widehat{\beta}^{(t)}+ \mathbf{v}(\widehat{\beta}^{(t)})dt$$
\item Statement Two (statistical description of fluid): Huge amounts of $ \widehat{\beta}^{(\cdot)} $'s move simultaneously and construct a statistical density $ \rho(x;t) $ at time $ t $ and location $ x $. $ \rho(x;t) $ satisfies $$
\frac{\partial}{\partial t}\rho(x;t) + \mathrm{div}\left(\mathbf{v}(x)\rho(x;t)\right)=0
$$
\end{itemize}

\begin{lemma}
Statement One is equivalent with Statement Two.
\end{lemma}

\begin{prof}\label{WrongProof}
We start by proving Statement One contains Statement Two. Consider location domain $ \Omega $ and time period $ (t,t+dt] $. $ \int_\Omega \frac{\partial\rho}{\partial t} dx\cdot dt$ is the amount of charge flowing through in $ \Omega $ during time $ (t,t+dt] $. Let $ \Omega_1 $ be the spatial location where particles flowing in are distributed at time $ t+dt $, and $ \Omega_2 $ be the spatial location where particles flowing out are distributed at time $ t+dt $. $ dt $ is an extremely short period such that particles either flow in or flow out. $ \Omega_1 $ and $ \Omega_2 $ do not intercept, and $ \Omega_1 + \Omega_2 $ constructs an extremely narrow boundary of the considered location domain $ \Omega $. We have
\begin{align*}
&\int_\Omega \frac{\partial\rho}{\partial t} dx\cdot dt\\
=&\int_{\Omega_1}\rho(x;t+dt)dx - \int_{\Omega_2}\rho(x;t+dt)dx
\end{align*}
\includegraphics[width=0.95\textwidth]{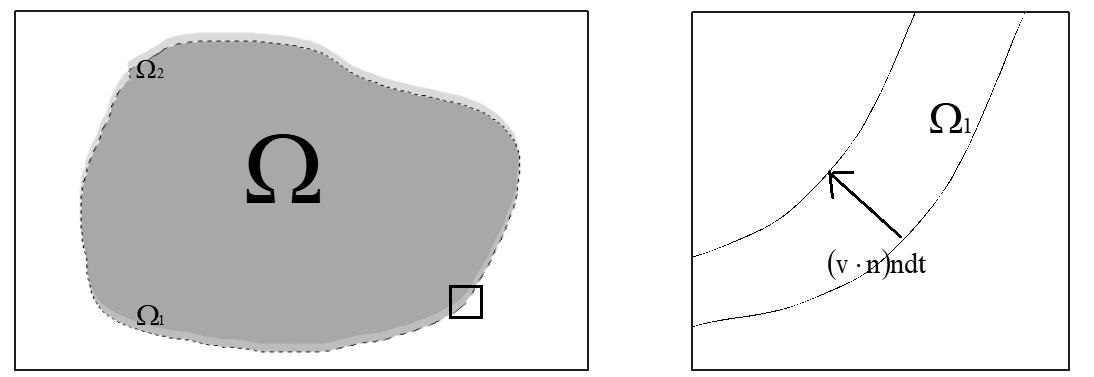}\\
Because $ dt $ is extremely short, $ \Omega_1 $ and $ \Omega_2 $ are very narrow boundary, any locality of which can be seen as rectangular infinitesimal element. We proceed
\begin{align}
&\int_{\Omega_1}\rho(x;t+dt)dx - \int_{\Omega_2}\rho(x;t+dt)dx \notag\\
=&\int_{\partial\Omega\cap\Omega_1} \left( \int_{y\in[s,s+(\mathbf{v}(s)\cdot\mathbf{n})\mathbf{n}dt]} \rho(y;t+dt)dy \right)dS \notag\\ 
&- \int_{\partial\Omega\cap\Omega_2} \left( \int_{y\in[s,s+(\mathbf{v}(s)\cdot\mathbf{n})\mathbf{n}dt]} \rho(y;t+dt)dy \right)dS \label{5y2}
\end{align}
Here $ s $ is a point on $ \partial\Omega $, $ [s,s+(\mathbf{v}(s)\cdot\mathbf{n})\mathbf{n}dt] $ indicates variable $ y $ taking integration over this one dimensional line segment, and $ \mathbf{n} $ is the outward unit normal vector of the surface $ \partial\Omega $. The readers may physically regard the progress as integrating $ \rho(\cdot;t+dt) $ over a surface $ \partial\Omega $ which is of little thickness, but mathematically it happens substitution of variable. A serious problem arises that this substitution of variable is nowhere bijection if the surface curves. It can be fixed in two ways: One is to consider rectangular domains $ \Omega $ only, and the measure theory later holds if only rectangles are concerned. The other is to strengthen constraint of the curvature of $ \partial\Omega $, the boundary of $ \Omega $, and yields a $ o(dt) $.

As $$
\rho(y,t+dt) = \rho(x,t+dt) + \triangledown_x\rho(x+\xi (y-x),t+dt)\vert_{\xi\in[0,1]}\cdot(y-x)
$$and by integrating $ \triangledown_x\rho(x+\xi (y-x),t+dt)(y-x) $ yields
\begin{align*}
&\int_{y\in[x,x+(\mathbf{v}(x)\cdot\mathbf{n})\mathbf{n}dt]} \triangledown_x\rho(x+\xi (y-x),t+dt)\cdot(y-x) dy \\
\leqslant& M \int_{y\in[x,x+(\mathbf{v}(x)\cdot\mathbf{n})\mathbf{n}dt]}(y-x) dy \\
=& M(dt)^2
\end{align*}where $ M $ is an absorbing constant, here and in all text we omit writing the hidden $ C^1 $ pre-condition of $ \rho(x,t) $. We proceed from \eqref{5y2} to:
\begin{align*}
=&\int_{\partial\Omega\cap\Omega_1} \left( \int_{y\in[x,x+(\mathbf{v}(x)\cdot\mathbf{n})\mathbf{n}dt]} \rho(x;t+dt)dy +o(dt)\right)dS \\
&- \int_{\partial\Omega\cap\Omega_2} \left( \int_{y\in[x,x+(\mathbf{v}(x)\cdot\mathbf{n})\mathbf{n}dt]} \rho(x;t+dt)dy +o(dt)\right)dS\\
=& \int_{\partial\Omega\cap\Omega_1 + \partial\Omega\cap\Omega_2} - \mathbf{n}\cdot \mathbf{v}(x)\rho(x;t+dt)dSdt+o(dt)\\
=& \int_{\partial\Omega} - \mathbf{n}\cdot \mathbf{v}(x)\rho(x;t+dt)dSdt+o(dt)\\
=&\int_\Omega -\mathrm{div}( \mathbf{v}(x)\rho(x;t+dt))dxdt+o(dt)
\end{align*}
Notice at $ \partial\Omega\cap\Omega_1 $, $ \mathbf{v}(x)\cdot\mathbf{n} $ is negative and there is additional negative sign yielded. Now we have proven that \begin{equation}\label{eq4.1.1}
\int_\Omega \frac{\partial\rho}{\partial t} dx dt = \int_\Omega -\mathrm{div}( \mathbf{v}(x)\rho(x;t+dt))dxdt+o(dt)
\end{equation}
by integrating \eqref{eq4.1.1} from time $ t_1 $ to time $ t_2 $ $$
\int_{[t_1,t_2]}\int_\Omega \frac{\partial\rho}{\partial t} dx dt = \int_{[t_1,t_2]}\int_\Omega -\mathrm{div}( \mathbf{v}(x)\rho(x;t+dt))dxdt
$$
At last, location-time domains like $ \Omega\times[t_1,t_2] $ constitutes semi-algebra in measure theory and thus $$
\iint_A \frac{\partial\rho}{\partial t} dx dt = \iint_A -\mathrm{div}( \mathbf{v}(x)\rho(x;t))dxdt
$$for any location-time $ A $, which is equivalent of $$
\frac{\partial\rho}{\partial t} + \mathrm{div}( \mathbf{v}(x)\rho(x;t)) = 0
$$This completes the proof of Statement One containing Statement Two, particle nature containing statistical nature.

In the contrary, by proceeding statistical nature contains particle nature, it is always attractive to consider Dirac delta functions. We want PDE to describe a scenario that huge amount of particles concentrating at one point. We advance the result by giving a testimonial but a failed proof in the following four steps:

\begin{itemize}
\item Step 1: Write PDE$$
\frac{\partial}{\partial t}\rho(x;t) + \mathrm{div}\left(\mathbf{v}(x)\rho(x;t)\right)=0
$$into its differentiation form:\begin{equation}\label{differentiationFormOfPDE}
\rho(x,t+dt) = \rho(x,t) - div(\mathbf{v}(x)\rho(x,t)) dt
\end{equation}

\item Step 2: Substitute $ \rho(x,t)= \delta_x $ into the right hand side of  \eqref{differentiationFormOfPDE}
\begin{align*}
\rho(x,t+dt) =& \rho(x,t) - div(\mathbf{v}(x))\rho(x,t) dt-
\mathbf{v}(x)\cdot \triangledown\rho(x,t)dt \\
=& (1-div(\mathbf{v}(x))dt)\delta_x -\triangledown\delta_x\cdot\mathbf{v}(x)dt
\end{align*}

\item Step 3: Apply the mean value theorem of Delta functions:
We assert that:$$
\delta_{x+\mathbf{v}(x)dt} = \delta_x -\triangledown\delta_x\cdot\mathbf{v}(x)dt
$$This is because for any test function $ f $, there exists\begin{align*}
&\int f(y)(\delta_x -\triangledown\delta_x  \cdot\mathbf{v}(x)dt) dy\\
=& f(x) - \int f(y)\triangledown\delta_x dy \cdot\mathbf{v}(x)dt \\
=& f(x) + \left(\int \triangledown f(y)\delta_x dy\right)\cdot\mathbf{v}(x)dt \\
=& f(x) + \triangledown f(x) \cdot\mathbf{v}(x)dt\\
=& f(x + \mathbf{v}(x)dt)\\
=& \int f(y)\delta_{x+\mathbf{v}(x)dt} dy
\end{align*}
Thus $$
\rho(x,t+dt) = \delta_{x+\mathbf{v}(x)dt} - div(\mathbf{v}(x))dt\delta_x
$$
\item Step 4: $ - div(\mathbf{v}(x))dt\delta_x = O(dt)$. Velocity field $ \mathbf{v} $ is physically smooth, meaning it does not vary dramatically in any locality. The gradient of $ \mathbf{v} $ is believed to be practically bounded.
\end{itemize}\end{prof}

The proof of Statement Two implying Statement One in the last lemma is not a mathematically acceptable proof but strictly speaking more an illustration. There are two deflects, namely:
\begin{itemize}
\item $ \rho(x,t+dt) $ is $ \delta_{x+\mathbf{v}(x)dt} $ but with an additional term $ O(dt) $
\item The mean value theorem of delta functions is applied in physics but remain imprecise in its mathematical theory in the author's knowledge.
\end{itemize}
Thus we have only proved that particle nature implying statistical nature, and we leave arduous proof of the contrary in the next section.

\section{Main Results}
The relationship between Newton-Raphson algorithms (or its variant Fisher-scoring algorithm in mathematical statistics) and the continuity equations in fluid dynamics (in its more well-known format the PDE of conservation of electric charge) is summarised in one sentence: the continuity equations is the statistical format of Newton's method, with the latter being the particle format of the former. They address the same mathematical phenomenon.

In this section and the text later, we write $ S = -F $ as the target function in Newton's method or the score function in mathematical statistics. We write $ I $ as the gradient of the target function $ F $ in Newton's method or the Fisher matrix in statistics. We write $ \widehat{\beta}^{(\cdot)} $ as the root solver in Newton's method or the maximum likelihood estimator in statistics. In many applications of Newton's method, many start points are tried to solve all the roots or to identify the global minimizer. This is a picture that many $ \widehat{\beta}^{(\cdot)} $'s move in the velocity field $ I(\cdot)^{-1}S(\cdot) $ independently and simultaneously, motivating us to consider the statistical scenario of many $ \widehat{\beta}^{(\cdot)} $'s.

By stating equivalence in mathematics words there is bijection between a Newton's method $$
\dot{\beta} = (I^{-1}S)(\beta)
$$ taking different start points, and a family of PDE $$
\frac{\partial}{\partial t}\rho(x;t) + \mathrm{div}\left(I(x)^{-1}S(x)\rho(x;t)\right)=0
$$ taking various initial value $ \rho(x;0) \in C^1$ each.

Sequential Newton's method \begin{equation}\label{FS_time}
\widehat{\beta}^{(t+1)}=\widehat{\beta}^{(t)}+ I(\widehat{\beta}^{(t)})^{-1}S(\widehat{\beta}^{(t)})
\end{equation} writes in its continuous-time form \begin{equation}\label{FS_cont_time}
\widehat{\beta}^{(t+dt)}=\widehat{\beta}^{(t)}+ I(\widehat{\beta}^{(t)})^{-1}S(\widehat{\beta}^{(t)})dt
\end{equation}
It is obviously seen that \eqref{FS_time} is the discrete-time numerical solution taking $ dt = 1 $ of \eqref{FS_cont_time}. The existence and uniqueness of ODE of \eqref{FS_cont_time} is characterised by a methodology of Pichard iterating sequence with its sufficient condition $ I^{-1}S $ being Lipschitz continuous. In this text continuous-time Newton's method \eqref{FS_cont_time} is what we concern.

\begin{theorem}\label{FS_is_PDE}
When the root solver or the maximum likelihood estimator $ \widehat{\beta}^{(\cdot)} $ runs in \eqref{FS_cont_time}, consider many $ \widehat{\beta}^{(t)} $'s move simultaneously. Huge amount of $ \widehat{\beta}^{(\cdot)} $ at time $ t $, $ \widehat{\beta}^{(t)} $ constitute macro statistical density $ \rho(x;t) $. Then $ \rho(x;t) $ satisfies PDE $$
\frac{\partial}{\partial t}\rho(x;t) + \mathrm{div}\left(I(x)^{-1}S(x)\rho(x;t)\right)=0
$$
\end{theorem}

\begin{prof}
By taking $ \mathbf{v}(x) = I(x)^{-1}S(x) $ in the lemma \ref{WrongProof} the theorem is proven.
\end{prof}

Theorem \ref{FS_is_PDE} justifies that the macro statistical format of the Newton-Raphson algorithm is the continuity equations. Namely when particles move in the velocity field\begin{equation}\label{FS2}
\widehat{\beta}^{(t+dt)}=\widehat{\beta}^{(t)}+ I(\widehat{\beta}^{(t)})^{-1}S(\widehat{\beta}^{(t)})dt
\end{equation} its statistical density $ \rho(x,t) $ satisfies \begin{equation}\label{PDE2}
\frac{\partial}{\partial t}\rho(x,t) + \mathrm{div}\left(I(x)^{-1}S(x)\rho(x,t)\right)=0
\end{equation}
We proceed proving that the statistical nature contains the particle nature. Since Lemma \ref{WrongProof} fails to be a complete proof, we proceed this phenomenon by the following one theorem and two lemmas.

\begin{theorem}\label{PDE_is_FS}
Particles $ \widehat{\beta}^{(t)} $'s move in certain velocity field $ \mathbf{v}(\cdot) $ with macro statistical density $ \rho(x,t) $. The latter is characterised by the PDE$$
\frac{\partial}{\partial t}\rho(x,t) + \mathrm{div}\left(I(x)^{-1}S(x)\rho(x,t)\right)=0
$$Then we have:$$
\mathbf{v}(\cdot) = I(\cdot)^{-1}S(\cdot)
$$
\end{theorem}

\begin{prof}
We know from Lemma \ref{WrongProof} that particles moving in the velocity field $ \mathbf{v}(\cdot) $ has the statistical density $ \rho(x,t) $ obeying $$
\frac{\partial}{\partial t}\rho(x,t) + \mathrm{div}\left(\mathbf{v}(x)\rho(x,t)\right)=0
$$
Given $ \rho(x,t) $ also obeys $$
\frac{\partial}{\partial t}\rho(x,t) + \mathrm{div}\left(I(x)^{-1}S(x)\rho(x,t)\right)=0
$$
We have\begin{equation}\label{eq5.4.1}
\mathrm{div}\left(\left(I(x)^{-1}S(x)-\mathbf{v}(x)\right)\rho(x,t)\right)\equiv 0 
\end{equation}
Notice by stating the regarded statistics nature of particles of the root solvers or the minimizers, we mean that in any location and time, with regrad to any statistical density $ \rho(x,t) $, \eqref{eq5.4.1} holds. Rewrite $ \rho(x,t) \equiv \rho(x_1,x_2,\cdots,x_p,t) $, $ p $ is the dimension of location, and take $ \rho(x,t) $ as specific form $ \rho(x,t) \equiv \rho_j(x_j) $, we derive \begin{align*}
&\mathrm{div}\left(I(x)^{-1}S(x)-\mathbf{v}(x)\right)\rho_j(x_j) \\
&+ \left(I(x)^{-1}S(x)-\mathbf{v}(x)\right)\cdot \triangledown \rho_j(x_j)\equiv 0 
\end{align*}
By taking $ \rho\equiv 1 $ in \eqref{eq5.4.1} $ \mathrm{div}\left(\left(I(x)^{-1}S(x)-\mathbf{v}(x)\right)\right)\equiv 0  $, this becomes $$
\left(I(x)^{-1}S(x)-\mathbf{v}(x)\right)\cdot \triangledown \rho_j(x_j)\equiv 0 
$$which is \begin{equation}\label{eq5.4.2}
\left(I(x)^{-1}S(x)-\mathbf{v}(x)\right)_j \frac{d}{dx_j}\rho_j(x_j)\equiv 0
\end{equation}
where $ \left(I(x)^{-1}S(x)-\mathbf{v}(x)\right)_j $ denotes the $ j $-th component of $ \left(I(x)^{-1}S(x)-\mathbf{v}(x)\right) $. By the arbitrariness of $ \frac{d}{dx_j}\rho_j(x_j) $ in \eqref{eq5.4.2}, the proof of the theorem is completed.
\end{prof}

Theorem \ref{PDE_is_FS} justifies the uniqueness of a velocity field given statistical description \eqref{PDE2}. To further proceed completing the statistics nature containing particle nature, we need to prove that macro statistical PDE designates a velocity field beforehand. Borrowing words from probability, we want to prove that if\begin{equation}\label{PDE3}
\frac{\partial}{\partial t}\rho(x,t) + \mathrm{div}\left(\mathbf{v}(x)\rho(x,t)\right)=0
\end{equation}
holds for any macro density $ \rho(x,t) $ then the statement \begin{equation}\label{statement3}
\textbf{Every particle moves in the velocity field }\mathbf{v}(\cdot)\textbf{ with probability one.}
\end{equation}
holds.

Assume the expectation of $ \rho(\cdot;t) $ is $ x $ $$
\int u\rho(u;t)du = x
$$and the variance of $ \rho(\cdot;t) $ is $ \Sigma $ $$
\int (u-x)(u-x)^T\rho(u;t)du = \Sigma
$$Presume $ \rho(\cdot;t) $ is non-degenerate that $ \Sigma $ is strictly positive. Consider the scenario that $ \Sigma $ being very small $ \Sigma\approx 0 $. In this case $ \rho(\cdot;t) $ is almost a Dirac delta function which concentrates at $ x $. In the text that follows, $ o(\Sigma) $ means how small it is compared with $ \Sigma $ but not necessarily in certain limiting process. The readers may regard it happening a limiting process if sequentially consider many more and more shrinking $ \rho(\cdot;t) $ to the delta function, but in order to most clearly present the text we do not arduously write all these out.

\begin{lemma}\label{lem5.1}
The expectation of $ \rho(\cdot;t+dt) $ is $$
x + \mathbf{v}(x)dt + o(\Sigma)dt + o(dt)
$$
\end{lemma}
\begin{prof}
Rewrite \eqref{PDE3} in the differentiation form:\begin{equation}\label{PDE4}
\rho(u;t+dt) = \rho(u;t) - \mathrm{div}\left(\mathbf{v}(u)\rho(u,t)dt\right) + o(dt)
\end{equation}and take expectation of both sides of \eqref{PDE4}. $
\int u\rho(u;t)du = x
$ is given.$$
-\int u \mathrm{div}\left(\mathbf{v}(u)\rho(u,t)dt\right)du = \int \left(\mathbf{v}(u)dt\right)\rho(u,t)du
$$This is because $\triangledown u $ yields identity matrix. Notice here that the integrable variable is $ du $, $ dt $ is the differentiation already done before. Assert $$
\int \left(\mathbf{v}(u)dt\right)\rho(u,t)du = \mathbf{v}(x)dt + o(\Sigma)dt
$$and once we prove this assertion the lemma is proved. \begin{align}\label{eqn3}
& \int \left(\mathbf{v}(u)dt\right)\rho(u,t)du - \mathbf{v}(x)dt=\\ \notag
& \int_{B(x,\epsilon)} \left(\mathbf{v}(u)dt - \mathbf{v}(x)dt\right)\rho(u,t)du +\int_{B(x,\epsilon)^c}  \left(\mathbf{v}(u)dt- \mathbf{v}(x)dt\right)\rho(u,t)du 
\end{align}
We look back at the nature of the problem this time. Velocity field $ \mathbf{v} $ is given and any intensively shrinking density is on the stage. In the whole text we omit writing $ C^1 $ condition of $ \mathbf{v} $, because it is physic common sense. But physics respect more than $ C^1 $ of the velocity field. By indicating a field we mean the continuity and uniformity of it, namely we do not see dramatic vibration in any locality. For small $ \epsilon $ $ \mathbf{v} $ is unvarying and the first term is arbitrarily close to $ 0 $. Given that $ \Sigma\approx 0 $ and as in real numerical study we take the whole measure of $ \rho $ in a compact set, which is a stronger condition than the tightness of probability measure guaranteed, and by doing variable substitution $ u = x + \Sigma^{\frac{1}{2}}v $ it is clear to see the second term is always a $ o(\Sigma)dt $. Now we rewrite: $$
\int \left(\mathbf{v}(u)dt\right)\rho(u,t)du - \mathbf{v}(x)dt= \left(o(\epsilon) + o(\Sigma(\epsilon) )\right)dt
$$
The analysis above can be done in a different approach. Let us reconsider \eqref{eqn3}. Do Taylor expansion of the first term to second order and \eqref{eqn3} becomes\begin{align*}
&\int_{B(x,\epsilon)} \left(\triangledown \mathbf{v}(x)dt\right)\cdot(u-x)\rho(u,t)du \\+
&\int_{B(x,\epsilon)} \frac{1}{2}(u-x)^T H(x)(u-x)dt\rho(u,t)du \\+ 
&\int_{B(x,\epsilon)^c}  \left(\mathbf{v}(u)dt- \mathbf{v}(x)dt\right)\rho(u,t)du
\end{align*}
Let $ \lambda $ be the maximum eigenvalue of $ H $ in $ B(x,\epsilon) $ the second term$$
\int_{B(x,\epsilon)} (u-x)^T H(x)(u-x)\rho(u,t)du \leqslant \lambda\int_{B(x,\epsilon)} (u-x)^T(u-x)\rho(u,t)du = o(\epsilon)
$$Still by doing variable substitution $ u = x + \Sigma^{\frac{1}{2}}v $ the third term $$ \int_{B(x,\epsilon)^c}  \left(\mathbf{v}(u)dt- \mathbf{v}(x)dt\right)\rho(u,t)du $$ is $ o(\Sigma(\epsilon))dt $. Here there is a trick. We further absorb a$$
-\left(\triangledown \mathbf{v}(x)dt\right)\cdot\int_{B(x,\epsilon)^c} (u-x)\rho(u,t)du
$$ into $ o(\Sigma(\epsilon))dt $ the first term becomes $$
\left(\triangledown \mathbf{v}(x)dt\right)\cdot\int (u-x)\rho(u,t)du =0
$$\eqref{eqn3} is now $ \left(o(\epsilon) + o(\Sigma(\epsilon) )\right)dt $.

Both two approaches reveal the result $ \left(o(\epsilon) + o(\Sigma(\epsilon) )\right)dt dt$. Review the nature of the problem, arbitrarily small $ \Sigma $ is considered, because we want cases alike a Dirac delta function. Thus$$
o(\epsilon) + o(\Sigma(\epsilon) ) = o(\Sigma)
$$This completes the proof of the lemma.
\end{prof}

\begin{lemma}\label{lem5.2}
The variance of $ \rho(\cdot;t+dt) $ is $$
\Sigma + o(\Sigma)dt + o(dt)
$$
\end{lemma}
\begin{prof}
We proceed by analysing\begin{equation}\label{lem5.2varicance}
\int (u-x -  \mathbf{v}(x)dt - o(\Sigma)dt - o(dt))(u-x -  \mathbf{v}(x)dt - o(\Sigma)dt - o(dt))^T\rho(u;t+dt)du
\end{equation}still rewrite $ \rho(u;t+dt) $ in the differentiation form \eqref{PDE4} $$
\rho(u;t+dt) = \rho(u;t) - \mathrm{div}\left(\mathbf{v}(u)\rho(u,t)dt\right) + o(dt)
$$This is a little disastrous in writing but the analysis is clear:
\begin{itemize}
\item $ \int (u-x)(u-x)^T\rho(u;t)du $ yields $ \Sigma $.
\item A $ dt $ item $ \int (u-x)\mathbf{v}(x)^T\rho(u;t)du $ yields $ 0 $.
\item We need only to tackle the remaining $ dt $ term being $ -\int (u-x)(u-x)^T\mathrm{div}\left(\mathbf{v}(u)\rho(u,t)dt\right)du $ because the rest are all $ o(dt) $.
\item Consider the $ (i,j) $-th component of $ (u-x)(u-x)^T $ and consider its integration $ -\int (u-x)_i(u-x)_j\mathrm{div}\left(\mathbf{v}(u)\rho(u,t)dt\right)du $, which equals to$$
\int \left((u-x)_i \mathbf{v}(u)_j+ (u-x)_j \mathbf{v}(u)_i\right)(\rho(u,t)dt)du
$$Apply the same technique in Lemma \ref{lem5.2}:\begin{align*}
&\int (u-x)_i \mathbf{v}(u)_j\rho(u,t)du = \\
&\int_{B(x,\epsilon)} (u-x)_i \mathbf{v}(u)_j\rho(u,t)du+\int_{B(x,\epsilon)^c} (u-x)_i \mathbf{v}(u)_j\rho(u,t)du
\end{align*} and $ \int_{B(x,\epsilon)^c} (u-x)_i \mathbf{v}(u)_j\rho(u,t)du $ always being $ o(\Sigma(\epsilon)) $. By the arbitrariness of $ \epsilon $ and the smoothness precondition of the velocity field $$
\int_{B(x,\epsilon)} (u-x)_i \mathbf{v}(u)_j\rho(u,t)du = \int_{B(x,\epsilon)} (u-x)_i \mathbf{v}(x)_j\rho(u,t)du + o(\epsilon)
$$ By further provide a $ \int_{B(x,\epsilon)^c} (u-x)_i \mathbf{v}(x)_j\rho(u,t)du $ that can be absorbed into $ o(\Sigma(\epsilon)) $ the first term becomes $ 0 $, and finally$$
\int \left((u-x)_i \mathbf{v}(u)_j+ (u-x)_j \mathbf{v}(u)_i\right)(\rho(u,t))du =o(\epsilon) + o(\Sigma(\epsilon))= o(\Sigma)
$$ In fact by applying Chebyshev's inequality directly the $ \int (u-x)_i \mathbf{v}(u)_j\rho(u,t)du $ can be proved a $ O(\Sigma) $, here we use $ \Sigma $ being arbitrarily small matrix to pursuit $ o(\Sigma) $ result.
\end{itemize}
The lemma is proved.
\end{prof}

\section{Numerical Study}
The essence of this section is to showcase that our equivalent PDE methodology holds numerically. To accomplish this goal we simulate a logistic generalised linear model with true $ \beta^*=(-0.2, 0.2, -0.2)^T $. We generate $ n=200 $ samples of $ \{(x_{j1},x_{j2},x_{j3},y_j)\}_{j=1}^{200} $, satisfying$$
y_j = \left\{\begin{array}{lr}
1,~\text{with probability}~\frac{\exp(x_j^T\beta^*)}{1+\exp(x_j^T\beta^*)}\\
0,~\text{with probability}~\frac{1}{1+\exp(x_j^T\beta^*)}
\end{array} \right.
$$which proceeds for the goal of maximising the likelihood$$
L(\beta) = \prod_j \left(\frac{\exp(x_j^T\beta)}{1+\exp(x_j^T\beta)}\right)^{y_j}\left(\frac{1}{1+\exp(x_j^T\beta)}\right)^{1-y_j}
$$or the equivalent log-likelihood$$
l(\beta) = \sum_j\left[ y_jx_j^T\beta - \log\left(1+\exp(x_j^T\beta)\right) \right]
$$The Score, corresponding $ -F $ in Newton's method$$
S(\beta) =  \frac{\partial l }{\partial\beta} = \sum_j \left[y_j-\frac{\exp(x_j^T\beta)}{1+\exp(x_j^T\beta)}\right]x_j
$$Here and in all canonical link cases Fisher matrix is gradient matrix with a negative sign:$$
I(\beta) =-\frac{\partial S }{\partial\beta^T} = \sum_j \left(\frac{\exp(x_j^T\beta)}{(1+\exp(x_j^T\beta))^2}\right)x_jx_j^T
$$Consider iteration$$
\dot{\beta} = I^{-1}(\beta)S(\beta)
$$in its equivalent PDE form$$
\rho(\cdot;t+dt) = \rho(\cdot;t) - \mathrm{div}(I^{-1}(\cdot)S(\cdot)\rho(\cdot;t)) dt
$$A Gaussian $ \rho(\cdot;0) $ is chosen, and the numerical grid length is set to be $ dx = 1/20 = 0.05 $. $ dt $ is set $ 0.05 $.

We draw in each dimension of $ \beta^* $ a 2-D intercepted sectional picture at time 0, 0.5, 1.0, and 2.0. At time 0 it is perfectly Gaussian:

\includegraphics[width=0.95\textwidth]{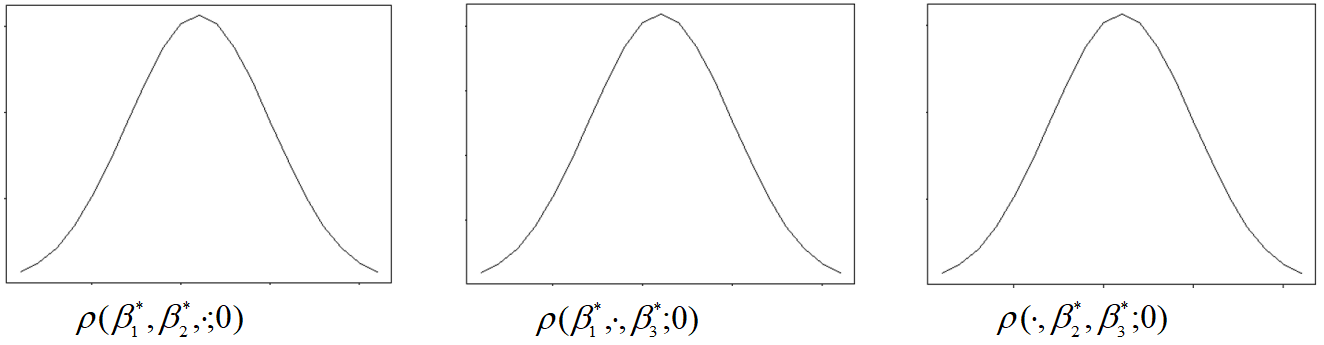}\\
At time 0.5 the convergence comes to light:

\includegraphics[width=0.95\textwidth]{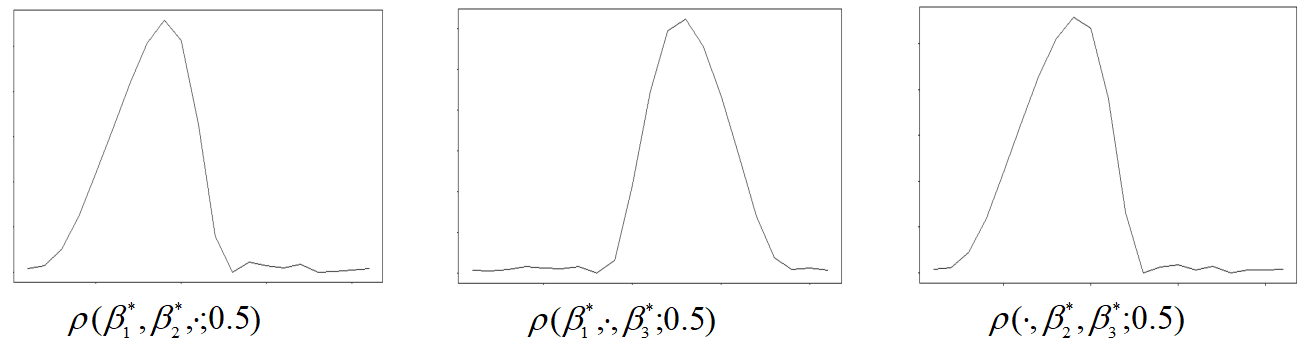}\\
At time 1 the convergence is clear:

\includegraphics[width=0.95\textwidth]{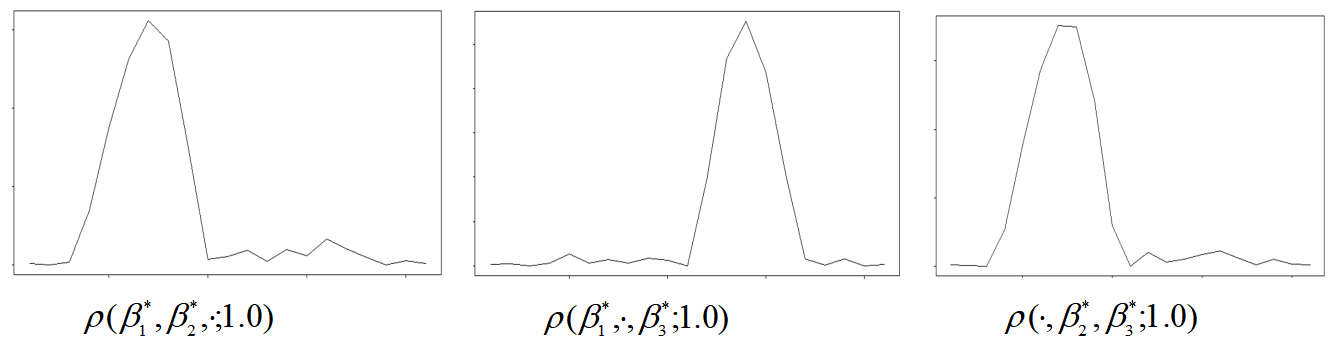}\\
At time 2 it is almost a Dirac delta function alike:

\includegraphics[width=0.95\textwidth]{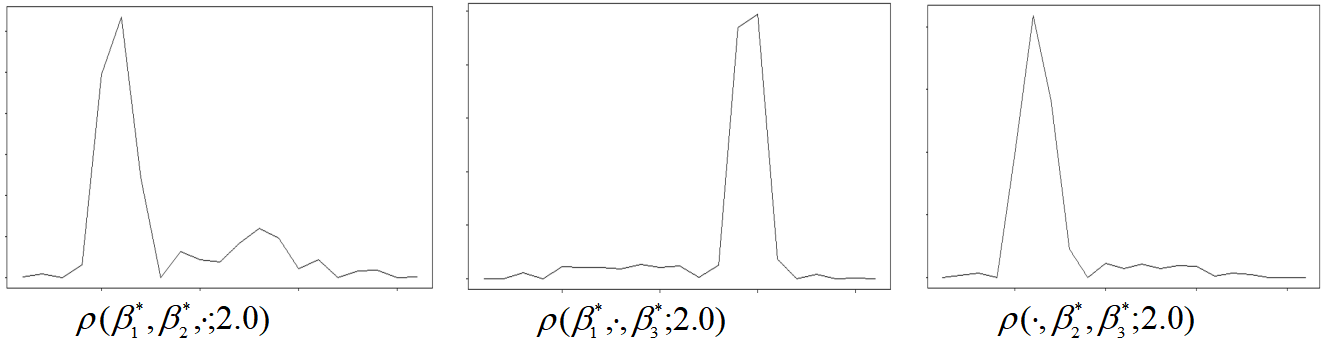}\\
To summarise, an overall movement picture is captured. Time consumption is about half a day. Curse of dimensionality is suffered, whereas Newton's method suffers from more sparsely distributed start points.

\section{Extensions}
May we use an $ itemize $ to structure this section. We propose the following four further studies:

\begin{itemize}
\item By stating equivalence we mean a bijection in mathematics. Particle motion in a velocity field is one-to-one correspondent to a family of continuity equations $$
\frac{\partial}{\partial t}\rho(x;t) + \mathrm{div}\left(\mathbf{v}(x)\rho(x;t)\right)=0
$$ taking miscellaneous initial $ \rho(x;0)\in C^1 $ each. Physically it is the equivalence of a macro velocity function with a micro velocity field. It follows direct questions of analytic, algebraic and geometric properties of the bijection.

\item Define total momenta of fluid as \begin{equation}
E(t)\triangleq \int \|\mathbf{v}(x)\|^2\rho(x,t)dx= \int \mathbf{v}(x)^T\mathbf{v}(x)\rho(x,t)dx
\end{equation}
with a $ \rho(x,t) $ obeying $$
\frac{\partial}{\partial t}\rho(x;t) + \mathrm{div}\left(\mathbf{v}(x)\rho(x;t)\right) = 0
$$Then we have\begin{align*}
\frac{d}{dt}E(t) =& \int \mathbf{v}(x)^T\mathbf{v}(x)\frac{\partial}{\partial t}\rho(x,t)dx\\
=&  - \int \mathbf{v}(x)^T\mathbf{v}(x)\mathrm{div}\left(\mathbf{v}(x)\rho(x;t)\right))dx\\
=& \int \triangledown\left(\mathbf{v}(x)^T\mathbf{v}(x)\right)\cdot\mathbf{v}(x)\rho(x;t)dx\\
=& 2 \int \mathbf{v}(x)^T\triangledown(\mathbf{v}(x))\mathbf{v}(x)\rho(x;t)dx
\end{align*}
When $ \triangledown(\mathbf{v}(x)) $ is negatively definite, total momenta of fluid decreases. There is exactly a case that $ \triangledown(\mathbf{v}(x)) $ is negatively definite, taking $ \mathbf{v}(x) =S(x)$, the score function in mathematical statistics. The result is characterised by Lemma \ref{Bartlett2} previously, though it is not capable of proving the convergence of the Fisher-scoring method.

\item Physics may be combined. The Maxwell equations dominates the three dimension. Two physic quantities electric field $ \mathbf{E} $ and magnetic field $ \mathbf{B} $ govern the whole system, and conservation of electric charge is encompassed.$$
\frac{\partial }{\partial t}\rho = \frac{\partial }{\partial t}\epsilon_0 \mathrm{div}\mathbf{E} =\epsilon_0 \mathrm{div}\frac{\partial }{\partial t}\mathbf{E}
$$By substituting $ \frac{\partial }{\partial t}\mathbf{E} = \frac{1}{\epsilon_0} \left( \mathrm{rot}\mathbf{B}\cdot\frac{1}{\mu_0} -\mathbf{j}\right)$ to the right hand side, this equals$$
\mathrm{div}\left( \mathrm{rot}\mathbf{B}\cdot\frac{1}{\mu_0} -\mathbf{j}\right) = -\mathrm{div}\mathbf{j}
$$Thus Newton's method is further equivalent with an electromagnetic system.

Quantum mechanics governs any arbitrary $ p $ dimension. In quantum mechanics $ \rho(x,t) $ is given by$$
\rho(x,t) = \psi^*(x,t)\psi(x,t)
$$Compare \eqref{continuity_equation} and \eqref{Schrodinger1} \begin{equation}\label{Schrodinger2}
\frac{\hbar}{i}\frac{1}{2m}\mathrm{div}\left[ \psi^*\triangledown \psi - \psi\triangledown \psi^*\right] = \mathrm{div}\left( \mathbf{v}(x)\psi^*\psi\right)
\end{equation}or$$
\mathrm{div}\left[ \psi^*\left(\frac{\hbar}{i}\frac{1}{m}\triangledown- \mathbf{v}(x)\right)\psi \right]= 
\mathrm{div}\left[ \psi\left(\frac{\hbar}{i}\frac{1}{m}\triangledown + \mathbf{v}(x)\right)\psi^*\right]
$$The right hand side equals to $ -\mathrm{div}\overline{\left[ \psi^*\left(\frac{\hbar}{i}\frac{1}{m}\triangledown - \mathbf{v}(x)\right)\psi\right]} $, and \eqref{Schrodinger2} becomes\begin{equation}\label{Schrodinger3}
\mathrm{div}~\mathrm{Re}\left[\psi^* \left(\frac{\widehat{p}}{m}-\mathbf{v}(x) \right)\psi\right]=0
\end{equation}
where $ \widehat{p} \equiv \frac{\hbar}{i}\triangledown $ is the momentum operator. Taking $ \mathbf{v}(x) = I^{-1}(x)S(x) $, \eqref{Schrodinger3} shapes quantum mechanics solution of Newton's method. The potential of the Schr\"{o}dinger equation of $ \psi $ needs to be real.

\item The continuous-time Fisher-scoring method \eqref{FS2} and \eqref{PDE2} may drop $ I^{-1} $ and consider $$
\widehat{\beta}^{(t+dt)}=\widehat{\beta}^{(t)}+ S(\widehat{\beta}^{(t)})dt
$$and its PDE counterpart$$
\frac{\partial}{\partial t}\rho(x,t) + \mathrm{div}\left(S(x)\rho(x,t)\right)=0
$$ Usually by considering a score we consider $ \frac{1}{n} $ times the score, it is just that $ \frac{1}{n} $ can be omitted from concern. Sometimes it is the root of the score to be concerned, and sometimes $ \frac{1}{n} $ is also putted before $ I $ and the two deduce. In many statistics studies, Fisher information matrix $ I $ is the gradient of $ S $ with negative sign, with an additional 0-expectation item, see Lemma \ref{Bartlett2}. $ I $ is positive definite in all cases. Thus for large $ n $, the gradient of $ S $ is negative definite. Local convergence of local maximum points is captured as long as numerical differentiation grid length being sufficiently small. To explain the phenomenon, let $ \beta_0 $ be a local maximum point and $ \beta_0 + \alpha\mathbf{l} $ be a position close to $ \beta_0 $, namely $ \alpha $ being small and $ \mathbf{l} $ unit directional vector. The velocity of $ \beta_0 + \alpha\mathbf{l} $ is $ S(\beta_0 + \alpha\mathbf{l}) $. Performing inner product with $ \mathbf{l} $, $$
\mathbf{l}^T S(\beta_0 + \alpha\mathbf{l}) = \mathbf{l}^T  \mathrm{grad} S(\beta_0) \mathbf{l} + o(\alpha) \leqslant 0 
$$and it can be seen that the velocity of $ S(\beta_0 + \alpha\mathbf{l}) $ is pointing inward towards $ \beta_0 $.
\end{itemize}

\end{document}